\documentclass[10pt]{article}

\usepackage[a4paper]{geometry}
\usepackage{indentfirst}
\usepackage{hyperref}
\usepackage{amssymb, amsmath, amsthm}

\hypersetup{bookmarksnumbered, pdfstartview=FitH, pdfpagemode=UseNone}

\newcommand{\llbr}{[\![}
\newcommand{\rrbr}{]\!]}

\newcommand{\bfu}{\mathbf{u}}
\newcommand{\bfv}{\mathbf{v}}
\newcommand{\bfw}{\mathbf{w}}
\newcommand{\bfx}{\mathbf{x}}
\newcommand{\bfy}{\mathbf{y}}

\newcommand{\sfW}{\mathsf{W}}
\newcommand{\sfB}{\mathsf{B}}
\newcommand{\sfM}{\mathsf{M}}

\newcommand{\sumW}{\mathbb{W}}
\newcommand{\sumB}{\mathbb{B}}
\newcommand{\sumM}{\mathbb{M}}
\newcommand{\sumU}{\mathbb{U}}

\newcommand{\bbD}{\mathbb{D}}
\newcommand{\bbE}{\mathbb{E}}
\newcommand{\bbF}{\mathbb{F}}

\newcommand{\bfalpha}{\boldsymbol\alpha}
\newcommand{\bfbeta}{\boldsymbol\beta}

\DeclareMathOperator{\score}{score}
\DeclareMathOperator{\surp}{surplus}

\newcommand{\mynewtheorem}[2]{\newtheorem{#1}{\indent #2}}
\newcommand{\myalttheorem}[2]{\newtheorem*{#1}{\indent #2}}
\mynewtheorem{conjecture}{Conjecture}
\myalttheorem{conjecture*}{Conjecture}
\mynewtheorem{lemma}{Lemma}
\mynewtheorem{proposition}{Proposition}
\myalttheorem{proposition*}{Proposition}
\mynewtheorem{corollary}{Corollary}
\myalttheorem{corollary*}{Corollary}
\mynewtheorem{theorem}{Theorem}
\myalttheorem{theorem*}{Theorem}
\newenvironment{myproof}[1][Proof]{\begin{proof}[\indent #1]}{\end{proof}}

\begin{document}

\title{\textbf{Powers of 2 in Balanced Grid Colourings}}
\author{Nikolai Beluhov}
\date{}

\maketitle

\begin{abstract} Let $B(m, n)$ be the number of ways to colour a $2m \times 2n$ grid in black and white so that, in each row and each column, half of the cells are white and half are black. Bhattacharya conjectured that the exponent of $2$ in the prime factorisation of $B(m, n)$ equals $s_2(m)s_2(n)$, where $s_2(x)$ denotes the number of $1$s in the binary expansion of $x$. We confirm this conjecture in some infinite families of special cases; most significantly, when $m$ is of the form either $2^k$ or $2^k + 1$ and $n$ is arbitrary. The proof when $m = 2^k + 1$ is substantially more difficult, and in connection with it we develop some general techniques for the analysis of inequalities between binary digit sums. \end{abstract}

\section{Introduction} \label{intro}

Let $m$ and $n$ be positive integers. Consider a rectangular grid of size $2m \times 2n$ whose cells are coloured in black and white. We call the colouring \emph{balanced} when, in each row and each column, half of the cells are white and half are black. Let $B(m, n)$ be the number of such colourings.

A different combinatorial interpretation of $B(m, n)$ is possible in graph theory. An \emph{Eulerian orientation} of a graph $G$ is one where, for each vertex, the in-degree equals the out-degree. Then $B(m, n)$ counts also the Eulerian orientations of the complete bipartite graph $K_{2m, 2n}$.

In the study of these numbers, it is often helpful to fix $m$ and consider the sequence obtained as $n$ varies. We denote this sequence by $\mathcal{B}(m)$. When $m = 1$, clearly $B(1, n) = \binom{2n}{n}$, and we get the central binomial coefficients. When $m = 2$, it is not too difficult to check that $B(2, n)$ equals the number of $2n$-step walks in the three-dimensional integer lattice which begin and end at the origin.~\cite{O1} The sequences we obtain when $m \ge 3$ are more mysterious. The On-Line Encyclopedia of Integer Sequences (henceforth OEIS) contains entries for the ones with $m = 3$, $4$, $\ldots$, $9$; as well as the entry \cite{O2} for the full two-parameter array $\{B(m, n)\}_{m, n \in \mathbb{N}}$ which links back to them.

The numbers $B(m, n)$ are discussed in \cite{N} and \cite{DBZ}. The two works use the graph-theoretic and the grid-based definitions, respectively. It is shown in \cite{DBZ} that each individual sequence $\mathcal{B}(m)$ is ``P-recursive''. Roughly speaking, P-recurrences are a generalisation of linear recurrences where the coefficients of the recurrence relation are not numerical constants but are allowed to depend polynomially on $n$ instead. However, the P-recurrences satisfied by the sequences $\mathcal{B}(m)$ seem to grow more complicated quite fast as $m$ increases; \cite{DBZ} gives them explicitly only when $m = 1$,~$2$,~$3$.

There are some classical enumeration problems where the exponent of $2$ in the prime factorisation of the answer exhibits sharp patterns which are difficult to verify in full generality. See, e.g., \cite{BL} for a discussion of this phenomenon in the context of domino tilings; or, in graph-theoretic language, the perfect matchings of grid graphs and their subgraphs.

Bhattacharya \cite{B} conjectured that the exponent of $2$ in $B(m, n)$ follows one remarkable pattern which has to do with the binary expansions of $m$ and $n$. Given a positive integer $x$, let $s_2(x)$ denote the binary digit sum of $x$; or, equivalently, the number of $1$s in the binary expansion of $x$.

\begin{conjecture*} \textnormal{(Bhattacharya)} The exponent of $2$ in the prime factorisation of the number of balanced grid colourings $B(m, n)$ equals $s_2(m)s_2(n)$. \end{conjecture*}

The numerical data available on OEIS agrees with the conjecture for all pairs $\{m, n\}$ such that $1 \le m \le 5$ and $1 \le n \le 50$. For the sequences $\mathcal{B}(m)$ with $m = 6$, $7$, $8$, $9$, OEIS lists the terms up to and including, respectively, $n = 39$, $30$, $25$, $22$. Direct computation shows that the conjecture holds true of these values as well.

But this is still a finite number of special cases. What about some infinite families? The case when $m = 1$ is immediate by Kummer's theorem -- see Section \ref{init}. Furthermore, the case when both of $m$ and $n$ are powers of $2$ was known to Bhattacharya at the time of \cite{B}. To this one, we add some other sparse families in Section \ref{bound}. However, most of all we will be interested in confirming the conjecture with $m$ fixed and $n$ arbitrary; i.e., over infinite families of the form~$\mathcal{B}(m)$. The case when $m$ is a power of $2$ is not too hard, and we give a short proof of the following in Section~\ref{init}:

\begin{theorem} \label{2k} Suppose that $m$ is of the form $2^k$ with $k$ nonnegative. Then the exponent of $2$ in the prime factorisation of the number of balanced grid colourings $B(m, n)$ equals $s_2(n)$. \end{theorem}

Our main result, on the other hand, is the case where $m$ is one greater than a power of $2$:

\begin{theorem} \label{2k1} Suppose that $m$ is of the form $2^k + 1$ with $k$ positive. Then the exponent of $2$ in the prime factorisation of the number of balanced grid colourings $B(m, n)$ equals $2s_2(n)$. \end{theorem}

This case is substantially more difficult; our proof occupies Sections \ref{three}--\ref{eq}. It hinges on certain inequalities between binary digit sums. In order to tackle them, we develop some general techniques for the analysis of a wider class of analogously-shaped problems.

The rest of the paper is structured as follows: Section \ref{init} gathers some initial observations. In Section \ref{bound}, we establish a number of partial results which ought to raise our confidence in the truth of the conjecture. For the proofs, temporarily we adopt a generating-function point of view -- but the rest of our work will be purely combinatorial.

Section \ref{three} resolves the case when $m = 3$ and $n$ is arbitrary; i.e., the sequence $\mathcal{B}(3)$. This argument forms the ``seed'' out of which our proof of the full Theorem \ref{2k1} will grow -- a process which begins in Section \ref{domino}. Section \ref{alg} is the aforementioned discussion of inequalities between binary digit sums. Then, in Sections \ref{ineq} and \ref{eq}, we apply what we have learned so as to complete the proof of Theorem \ref{2k1}. The concluding Section \ref{further} is supposed to list some potential directions for further research; however, what it does achieve in practice is mostly to remind us how little we know about the conjecture.

\section{Initial Observations} \label{init}

One result of particular importance to us will be Kummer's theorem for the prime factor $2$. This is usually stated as ``the exponent of $2$ in $\binom{x + y}{x, y}$ equals the number of carries when $x$ and $y$ are added together in binary''. But for our purposes it will be more convenient to express the number of carries in terms of the binary digit sums of $x$ and $y$, as $s_2(x) + s_2(y) - s_2(x + y)$.

The theorem generalises to multinomial coefficients, too. Then it states that the exponent of $2$ in $\binom{x_1 + x_2 + \cdots + x_k}{x_1, x_2, \ldots, x_k}$ equals $s_2(x_1) + s_2(x_2) + \cdots + s_2(x_k) - s_2(x_1 + x_2 + \cdots + x_k)$.

We return now to balanced grid colourings. For concreteness, we assume that in an $a \times b$ grid the width is $a$ and the height is $b$; i.e., each row consists of $a$ cells and each column consists of $b$ cells.

Observe, to begin with, that the problem has a lot of natural automorphisms. For example, swapping the two colours in a balanced colouring produces a different balanced colouring. (In the graph-theoretic interpretation, this corresponds to reversing the orientations of all arcs.) Since this is an involution, the number of balanced colourings is always even.

We obtain some more helpful automorphisms by permuting the rows of the grid. (Imagine the grid sliced into horizontal strips of size $2m \times 1$ each; these are then reordered and glued back together into a new grid. In the graph-theoretic interpretation, this corresponds to permuting the vertices of the $2n$-sized part of $K_{2m, 2n}$.) We can permute the columns as well. Or we can permute both the rows and the columns -- notice that a row permutation and a column permutation always commute.

Let $B^\dagger(m, n)$ denote the number of balanced colourings in which the left half of the top row consists of $m$ white cells and its right half consists of $m$ black cells. By permuting the columns, we easily get that $B(m, n) = \binom{2m}{m}B^\dagger(m, n)$. Similarly, let $B^\ddagger(m, n)$ denote the number of balanced colourings in which the top half of the leftmost column consists of $n$ white cells and its bottom half consists of $n$ black cells. Then $B(m, n) = \binom{2n}{n}B^\ddagger(m, n)$. By Kummer's theorem, it follows that the exponent of $2$ in $B(m, n)$ is at least $\max\{s_2(m), s_2(n)\}$.

We can get a slightly better lower bound by combining row and column permutations. Let $B^\mathsection(m, n)$ denote the number of balanced colourings which are counted in both of $B^\dagger(m, n)$ and $B^\ddagger(m, n)$. It is straightforward to see, by permuting both rows and columns, that $B(m, n) = \frac{1}{2}\binom{2m}{m}\binom{2n}{n}B^\mathsection(m, n)$. Hence, the exponent of $2$ in $B(m, n)$ is at least $s_2(m) + s_2(n) - 1$. This is a definite improvement -- but it is still a far cry from the conjectured value of $s_2(m)s_2(n)$.

We get the feeling that part of the difficulty lies in the desired exponent being quite large. So it makes sense first of all to try those special cases where the conjectured value is reasonably small. We reproduce now Bhattacharya's argument for the case when $s_2(m) = s_2(n) = 1$. Suppose that both of $m$ and $n$ are powers of $2$. Then, in order to confirm the conjecture, it suffices to calculate $B(m, n)$ modulo $4$.

Consider any balanced colouring with $k$ distinct kinds of rows which occur with multiplicities $x_1$, $x_2$, $\ldots$, $x_k$. (We say that two rows are of the same kind when they are identically coloured. Of course, $x_1 + x_2 + \cdots + x_k = 2n$.) Its orbit under the action of the row permutations is of size $\binom{2n}{x_1, x_2, \ldots, x_k}$. Since $k \ge 2$ and $s_2(n) = 1$, by Kummer's theorem this multinomial coefficient will be divisible by $4$ unless $k = 2$ and $x_1 = x_2 = n$.

We are left to investigate the balanced colourings of the latter form. The two kinds of rows in the colouring must form a complementary pair; i.e., swapping the two colours must interchange them. There are $\frac{1}{2}\binom{2m}{m}$ such pairs; and, as above, each one of them yields $\binom{2n}{n}$ balanced colourings. Since both of $m$ and $n$ are powers of $2$, by Kummer's theorem the former number is odd whereas the latter one works out to $2$ modulo $4$. This takes care of the case when $s_2(m) = s_2(n) = 1$.

We round off this section by resolving also the case when $s_2(m) = 1$ but no restrictions are imposed on $s_2(n)$ and $n$.

\begin{myproof}[Proof of Theorem \ref{2k}] Suppose that $m$ is a power of $2$. Since $B(m, n) = \binom{2n}{n} B^\ddagger(m, n)$ and the exponent of $2$ in $\binom{2n}{n}$ equals $s_2(n)$ by Kummer's theorem, what remains is to show that $B^\ddagger(m, n)$ is odd.

Consider any balanced colouring counted in $B^\ddagger(m, n)$. Trim off the leftmost column and partition the rest of the grid into two smaller sub-grids of size $(2m - 1) \times n$ each. Call the top one $G'$ and the bottom one $G''$.

Given a row or column coloured in black and white, we define its \emph{discrepancy} to be the number of white cells minus the number of black cells. So a colouring is balanced if and only if all of the discrepancies vanish. The discrepancy of each row of $G'$ is $-1$ and the discrepancy of each row of $G''$ is $1$. Let $\bfu' \in \mathbb{Z}^{2m - 1}$ be the vector formed by the column discrepancies of $G'$, listed in order from left to right. Define $\bfu''$ similarly for $G''$. Clearly, $\bfu' + \bfu'' = \mathbf{0}$.

Given any vector $\bfu \in \mathbb{Z}^{2m - 1}$, let $\Theta(\bfu)$ be the number of ways to colour $G'$ in black and white so that each row is of discrepancy $-1$ and the column discrepancies are given by $\bfu$. Swapping the two colours, we see that $\Theta(\bfu)$ is also the number of ways to colour the cells of $G''$ so that each row is of discrepancy $1$ and the column discrepancies are given by $-\bfu$. Hence, \[B^\ddagger(m, n) = \sum_\bfu \Theta(\bfu)^2 \equiv \sum_\bfu \Theta(\bfu) \pmod 2,\] where the last step follows because $x^2 \equiv x \pmod 2$ for all integers $x$.

However, $\sum_\bfu \Theta(\bfu)$ is just the number of ways to colour $G'$ in black and white so that each row is of discrepancy $-1$. Since there are $\binom{2m - 1}{m}$ options for each individual row, we conclude that there exist $\binom{2m - 1}{m}^n$ such colourings. This number is of the same parity as $\binom{2m - 1}{m}$, which in light of $m$ being a power of $2$ becomes odd by Kummer's theorem. \end{myproof}

\section{The Bounded Parts of the Conjecture} \label{bound}

We define the $h$-\emph{bounded part} of the conjecture to be the claim ``the conjecture is true of all $m$ and $n$ with $\max\{s_2(m), s_2(n)\} \le h$''. Notice that, for each concrete $h$, in order to confirm the $h$-bounded part of the conjecture it is enough to understand the numbers $B(m, n)$ modulo $2^{h^2 + 1}$. This makes the bounded parts substantially easier to get a grasp on than the full conjecture.

For our analysis of $B(m, n)$ modulo a fixed power of $2$, it will be convenient to view the problem in terms of generating functions. Let $S$ be the Laurent polynomial in $x_1$, $x_2$, $\ldots$, $x_{2m}$ obtained as the sum of all monomials $x_1^{e_1}x_2^{e_2} \cdots x_{2m}^{e_{2m}}$ with $|e_1| = |e_2| = \cdots = |e_{2m}| = 1$ and $e_1 + e_2 + \cdots + e_{2m} = 0$. Then $B(m, n)$ is the constant term of $S^{2n}$.

A sequence given by the constant terms of the powers of a fixed Laurent polynomial is known as a \emph{constant-term sequence}. Our Lemma \ref{ab} below is a corollary of one more general result on the arithmetic properties of these sequences -- Lemma 1.2 of \cite{MV}. We omit the details of the derivation in favour of a self-contained argument tailored specifically to the setting of balanced grid colourings.

\begin{lemma} \label{ab} Suppose that $n = 2^{a + b + 1}n' + n''$ with $n'' < 2^b$. Then \[B(m, n) \equiv B(m, 2^an')B(m, n'') \pmod{2^{a + 2}}.\] \end{lemma}

The parameters $a$, $b$, $n'$, $n''$ are all assumed to be nonnegative integers. Of course, $B(0, 0) = B(0, n) = B(m, 0) = 1$.

\begin{myproof} Let $P$ be any integer-coefficient Laurent polynomial in $x_1$, $x_2$, $\ldots$, $x_{2m}$. We write $P_{\langle e \rangle}$ for $P(x_1^e, x_2^e, \ldots, x_{2m}^e)$. By induction on $i$, it is straightforward to see that $P^{2^{i + 1}} = P_{\langle 2 \rangle}^{2^i} + 2^{i + 1}Q$ for some integer-coefficient Laurent polynomial $Q$ in $x_1$, $x_2$, $\ldots$, $x_{2m}$. Hence, \[S^{2^{a + b + 2}} \equiv S_{\langle 2 \rangle}^{2^{a + b + 1}} \equiv S_{\langle 4 \rangle}^{2^{a + b}} \equiv \cdots \equiv S_{\langle 2^{b + 1} \rangle}^{2^{a + 1}} \pmod{2^{a + 2}}.\]

Let $S' = S_{\langle 2^{b + 1} \rangle}^{2^{a + 1}n'}$ and $S'' = S^{2n''}$. Then $S^{2n} = S^{2^{a + b + 2}n'}S^{2n''} \equiv S'S'' \pmod{2^{a + 2}}$. We conclude that $B(m, n)$ is congruent modulo $2^{a + 2}$ to the constant term of $S'S''$. However, the exponents of $x_1$, $x_2$, $\ldots$, $x_{2m}$ in all monomials of $S'$ are divisible by $2^{b + 1}$; while the exponents of the same variables in the monomials of $S''$ cannot exceed $2n'' < 2^{b + 1}$ by absolute value. Thus the constant term of $S'S''$ coincides with the product of the constant terms of $S'$ and~$S''$. The former equals the constant term of $S^{2^{a + 1}n'}$, or $B(m, 2^an')$; whereas the latter equals $B(m, n'')$. \end{myproof}

Lemma \ref{ab} is most useful when the binary expansion of $n$ is ``sparse'', in the sense that the $1$s in it are few and far between. Let us formalise this notion. We say that a positive integer is $k$-\emph{sparse} when, for each pair of $1$s in its binary expansion, between these two $1$s there are at least $k$ zeroes.

\begin{proposition} \label{infn} For all $m$ and $h$, there exist infinitely many $n$ with $s_2(n) = h$ such that the conjecture is true of $m$ and $n$. \end{proposition}

\begin{myproof} Let $k = s_2(m) \cdot h$ and consider any $k$-sparse $n$ with $s_2(n) = h$. Suppose that $n = 2^{e_1} + 2^{e_2} + \cdots + 2^{e_h}$ with $e_1 > e_2 > \cdots > e_h$. By iterated application of Lemma \ref{ab}, we can split off these powers of $2$ one by one to find that \[B(m, n) \equiv B(m, 2^{e_1 - e_2 - 2})B(m, 2^{e_2 - e_3 - 2}) \cdots B(m, 2^{e_{h - 1} - e_h - 2})B(m, 2^{e_h}) \pmod{2^{k + 1}}.\]

However, the exponent of $2$ in the product on the right-hand side equals $s_2(m)s_2(n) < k + 1$ by Theorem \ref{2k}. \end{myproof}

In the next couple of propositions, ``suitable $m$ and $n$'' stands for ``$m$ and $n$ with $\max\{s_2(m),\allowbreak s_2(n)\} \le h$''; i.e., $m$ and $n$ which are relevant to the $h$-bounded part of the conjecture.

\begin{proposition} \label{short} Fix $h$ and suppose that the $h$-bounded part of the conjecture holds for all suitable $m$ and $n$ with at most $h^3$ binary digits. Then it holds for all suitable $m$ and $n$. \end{proposition}

So, for each fixed $h$, we can test whether the $h$-bounded part of the conjecture is true by means of a finite amount of computation. Though this ``finite amount'' is already prohibitively large when $h = 2$.

\begin{myproof} By induction on $m + n$. Suppose that $m$ and $n$ form a suitable pair with more than $h^3$ digits in the binary expansion of (for concreteness) $n$. Since there are at most $h$ instances of the digit $1$ and at most $h$ runs of the digit $0$ in this binary expansion, one of these runs (possibly the one at the end) must be of length at least $h^2$. Hence, $n = 2^{h^2 + b}n' + n''$ with $n' \neq 0$ and $n'' < 2^b$. By Lemma~\ref{ab}, it follows that \[B(m, n) \equiv B(m, 2^{h^2 - 1}n')B(m, n'') \pmod{2^{h^2 + 1}}.\]

What remains is to apply the induction hypothesis. \end{myproof}

We engage now in a bit of innocent sleight-of-hand.

\begin{proposition} \label{dense} For each fixed $h$, the $h$-bounded part of the conjecture is true of almost all suitable $m$ and $n$, in the following technical sense: The fraction of the suitable pairs $\{m, n\}$ for which it does hold, out of all suitable pairs $\{m, n\}$ with $\max\{m, n\} \le H$, approaches unity when $H$ grows without bound. \end{proposition}

We classify this result as ``sleight-of-hand'' because it sounds much grander than it really is. Yes, the $h$-bounded part of the conjecture does hold quite densely among the suitable pairs. But the suitable pairs themselves are, for each fixed $h$, spread out incredibly thinly among all pairs of positive integers.

\begin{myproof} It is a straightforward exercise to show that, as $H$ grows without bound, both of $m$ and $n$ are $h^2$-sparse in almost all suitable pairs $\{m, n\}$ with $\max\{m, n\} \le H$. However, for a pair which satisfies this condition, the $h$-bounded part of the conjecture holds by the same argument as in the proof of Proposition \ref{infn}. \end{myproof}

\section{The Case of \texorpdfstring{$m = 3$}{m = 3}} \label{three}

Throughout this section, let $m = 3$. The conjecture then specialises to the exponent of $2$ in $B(3, n)$ being equal to $2s_2(n)$.

Suppose we wish to show that the answer to some enumeration problem is divisible by a given positive integer $d$. One of the most natural approaches would be to look for a way to express that answer as a sum where each summand is divisible by $d$ on its own. This is just what we are going to do. There are two difficulties to overcome -- first, to come up with a suitable sum; and, second, to show that each summand in it is indeed divisible by $d$.

In the present setting, $d = 2^{2s_2(n)}$. Of course, in order to confirm the conjecture, divisibility by $d$ is not enough; we must also verify that $2s_2(n)$ is the exact value of the exponent. However, once we have achieved the former, the latter will turn out to follow as well without too much extra effort.

On to the proof. Partition the grid into horizontal dominoes; i.e., sub-grids of size $2 \times 1$. Each row will be split into three dominoes. Furthermore, the dominoes will form three ``domino columns'' each one of which is the union of two neighbouring columns of the original grid.

We label each domino with a letter out of the alphabet $\{\sfW, \sfB, \sfM\}$ depending on the colours of its cells, as follows: a pure white domino is labelled $\sfW$; a pure black domino is labelled $\sfB$; and a mixed domino, with two cells of different colours, is labelled $\sfM$. The \emph{type} of a row is the three-letter word formed by the labels of its dominoes, in order from left to right.

Since each row is balanced, the letters $\sfW$ and $\sfB$ must occur the same number of times in its type. So the feasible types in a balanced colouring are $\sfM\sfW\sfB$, $\sfM\sfB\sfW$, $\sfB\sfM\sfW$, $\sfW\sfM\sfB$, $\sfW\sfB\sfM$, $\sfB\sfW\sfM$, $\sfM\sfM\sfM$. Denote the number of rows of each type by $x_1$, $x_2$, $\ldots$, $x_7$, respectively. We call $\bfx = (x_1, x_2, \ldots, x_7)$ the \emph{multiplicity vector} of the colouring.

Now $x_1 + x_2 + \cdots + x_7 = 2n$. But this is not the only condition which $\bfx$ must satisfy. Since the colouring is balanced, the leftmost domino column must contain the same number of dominoes with the labels $\sfW$ and $\sfB$. The former are counted by $x_4 + x_5$ and the latter by $x_3 + x_6$. Reasoning analogously about the middle and rightmost domino columns, we arrive at the system \[\begin{gathered} x_1 + x_2 + \cdots + x_7 = 2n\\ x_1 + x_4 = x_2 + x_3\\ x_1 + x_6 = x_2 + x_5\\ \phantom{.}x_3 + x_6 = x_4 + x_5. \end{gathered} \tag{$\mathbf{A}_3$}\]

Let $\Omega$ be the set of all such $\bfx$. Given a multiplicity vector $\bfx \in \Omega$, how many balanced colourings are associated with it? The number of ways to assign types to rows is $\binom{2n}{x_1, x_2, \ldots, x_7}$. Furthermore, in each domino column, half of the mixed dominoes must go $\square\blacksquare$ and the other half must go $\blacksquare\square$. The three domino columns contain $x_1 + x_2 + x_7$, $x_3 + x_4 + x_7$, and $x_5 + x_6 + x_7$ mixed dominoes, respectively. We tally everything up, and arrive at a total of \[\binom{2n}{x_1, x_2, \ldots, x_7}\binom{x_1 + x_2 + x_7}{\frac{1}{2}(x_1 + x_2 + x_7)}\binom{x_3 + x_4 + x_7}{\frac{1}{2}(x_3 + x_4 + x_7)}\binom{x_5 + x_6 + x_7}{\frac{1}{2}(x_5 + x_6 + x_7)} \tag{$\mathbf{E}$}\] balanced colourings associated with $\bfx$.

We claim that the exponent of $2$ in ($\mathbf{E}$) is at least $2s_2(n)$ for all $\bfx \in \Omega$. This exponent is immediately expressible in terms of binary digit sums over $x_1$, $x_2$, $\ldots$, $x_7$ via Kummer's theorem. Our claim boils down to the inequality \[s_2(x_1) + s_2(x_2) + \cdots + s_2(x_7) + s_2(x_1 + x_2 + x_7) + s_2(x_3 + x_4 + x_7) + s_2(x_5 + x_6 + x_7) \ge 3s_2(n). \tag{$\mathbf{B}_3$}\]

Recall that $s_2(y) + s_2(z) \ge s_2(y + z)$ for all nonnegative integers $y$ and $z$. (E.g., because the difference between the two sides equals the number of carries when $y$ and $z$ are added together in binary -- as we remarked at the beginning of Section \ref{init}.) Taking ($\mathbf{A}_3$) into account, it follows~that \[s_2(x_1) + s_2(x_4) \ge s_2(x_1 + x_4) = s_2(2(x_1 + x_4)) = s_2(x_1 + x_2 + x_3 + x_4)\] and \begin{align*} s_2(x_1) + s_2(x_4) + s_2(x_5 + x_6 + x_7) &\ge s_2(x_1 + x_2 + x_3 + x_4) + s_2(x_5 + x_6 + x_7)\\ &\ge s_2(x_1 + x_2 + \cdots + x_7) = s_2(2n) = s_2(n). \end{align*}

Similarly, \begin{align*} s_2(x_2) + s_2(x_5) + s_2(x_3 + x_4 + x_7) &\ge s_2(n)\\ s_2(x_3) + s_2(x_6) + s_2(x_1 + x_2 + x_7) &\ge s_2(n). \end{align*}

Summing the latter three inequalities together with the obvious $s_2(x_7) \ge 0$, we obtain~($\mathbf{B}_3$). Our claim has been verified.

However, this establishes only a lower bound on the exponent of $2$ in $B(3, n)$. We must also show that our lower bound equals the exact value. To begin with, we can discard all $\bfx$ for which ($\mathbf{B}_3$) is strict. Call $\bfx$ \emph{exact} if it attains equality in ($\mathbf{B}_3$). Our goal will be achieved if we demonstrate that the number of exact multiplicity vectors is odd.

The colour-swapping automorphism of Section \ref{init} comes to our help. When we swap the two colours, the mixed dominoes remain mixed whereas the pure ones exchange labels; i.e., colour swapping induces over $\Omega$ the involution $(x_1, x_2, x_3, x_4, x_5, x_6, x_7) \leftrightarrow (x_2, x_1, x_4, x_3, x_6, x_5, x_7)$. Two multiplicity vectors exchanged by this involution always yield the same ($\mathbf{E}$). Hence, if one of them is exact, so will be the other one as well. We are left to analyse the exact multiplicity vectors which are preserved by the involution. Call them \emph{exceptional}.

Clearly, $\bfx$ is a fixed point of our involution if and only if $x_1 = x_2$, $x_3 = x_4$, $x_5 = x_6$. Denote these quantities by $a$, $b$, $c$, respectively. Our proof of ($\mathbf{B}_3$) implies that, for equality to be attained, $s_2(x_1) + s_2(x_4) = s_2(x_1 + x_4)$ must hold; i.e., there must be no carries when $a$ and $b$ are added together in binary. Similarly, no carries can occur in the binary evaluations of the sums $a + c$ and $b + c$, either. The attainment of equality in ($\mathbf{B}_3$) forces also $s_2(x_7) = 0$; or, equivalently, $x_7 = 0$.

Conversely, if $a$, $b$, $c$ are three nonnegative integers such that $a + b + c = n$ and the evaluation of the sum $a + b + c$ is carry-free in binary, then $(a, a, b, b, c, c, 0)$ will be an exceptional multiplicity vector. (The absence of carries in the three pairwise binary summations is equivalent to the absence of carries in the binary summation of the three quantities all together.) Or, in other words, the ordered triples $(a, b, c)$ which satisfy these conditions are in bijective correspondence with the exceptional multiplicity vectors. But how many such triples are there?

The binary expansion of $n$ reflects its partitioning into pairwise distinct powers of $2$. For the evaluation of the sum $a + b + c$ to be carry-free in binary, these powers of $2$ must be distributed between $a$, $b$, $c$. There are $s_2(n)$ items to distribute, and three ``bins'' to distribute them between. Thus the number of feasible triples becomes $3^{s_2(n)}$. This number is odd, and the proof of the conjecture in the special case when $m = 3$ is complete.

\section{Dominoes} \label{domino}

What changes when we attempt to run the argument of Section \ref{three} with $m \neq 3$? We can still partition the grid into horizontal dominoes. Each row will be split into $m$ of them, and the dominoes will also form $m$ domino columns. We define the type of a row as before; it will always be an $m$-letter word over the alphabet $\{\sfW, \sfB, \sfM\}$ in which the letters $\sfW$ and $\sfB$ occur equally many times.

Given a type $T$, let $x_T$ be the number of rows of type $T$ in our colouring. Fix an ordering of the types, and let $\bfx$ be the vector whose components are the multiplicities of the types, listed as in the ordering. This is how the notion of a multiplicity vector generalises.

We write $T[i]$ for the $i$-th letter of $T$. Let $\sumW_i$ be the sum of all multiplicities $x_T$ such that $T[i] = \sfW$. The sums $\sumB_i$ and $\sumM_i$ are defined similarly. These sums count the pure white, pure black, and mixed dominoes in the $i$-th domino column, respectively.

We also let $\sumU$ be the sum of $x_T$ over all $T$. Of course, we already know the numerical value of~$\sumU$ -- it must equal $2n$. So the reader might be puzzled as to why we are bothering to introduce special notation for it. The reason will become obvious in Section \ref{alg}. For the moment, note that we view both the preceding sums and $\sumU$ as formal expressions in the multiplicity variables rather than as numerical quantities.

Taking stock over each separate domino column, we find that $\bfx$ must satisfy the conditions \[\begin{gathered} \sumU = 2n\\ \sumW_i = \sumB_i \quad \text{for all $i$}. \end{gathered} \tag{$\mathbf{A}_m$}\]

Once again, we denote the set of all such $\bfx$ by $\Omega$. The calculation of the number of balanced colourings associated with a given multiplicity vector $\bfx \in \Omega$ generalises smoothly, too. We get one ``big'' multinomial coefficient multiplied by the product of $m$ ``smaller'' central binomial coefficients which govern the disambiguation of the mixed dominoes in the $m$ separate domino columns. By Kummer's theorem, the exponent of $2$ in this total works out to \[\sum_T s_2(x_T) - s_2(n) + \sum_i s_2(\sumM_i). \tag{$\mathbf{F}$}\]

How large can we guarantee this exponent to be? If $m$ is even, then the answer is ``not too large''. Consider the case where the only types which occur with nonzero multiplicities are $\sfW\sfW \ldots \sfW\sfB\sfB \ldots \sfB$ and $\sfB\sfB \ldots \sfB\sfW\sfW \ldots \sfW$. Naturally, both of them must occur with multiplicity~$n$. Then ($\mathbf{F}$) simplifies down to $s_2(n)$. But in Section \ref{init} we already got that the exponent of $2$ in $B(m, n)$ is at least $s_2(n)$ anyway -- so the ``domino method'' does not achieve much when $m$ is even.

What about the odd $m \ge 3$? Consider now the case where the only types which occur with nonzero multiplicities are $\sfW\sfW \ldots \sfW\sfM\sfB\sfB \ldots \sfB$ and $\sfB\sfB \ldots \sfB\sfM\sfW\sfW \ldots \sfW$. Once again, both of them must occur with multiplicity $n$. This time around, ($\mathbf{F}$) simplifies down to the more satisfactory $2s_2(n)$. We can never hope to show that ($\mathbf{F}$) is greater than $2s_2(n)$ in general; but, on the bright side, $2s_2(n)$ is still pretty good.

There remains the matter of finding out whether ($\mathbf{F}$) is in fact always at least $2s_2(n)$. Our desired inequality becomes \[\sum_T s_2(x_T) + \sum_i s_2(\sumM_i) \ge 3s_2(\sumU), \tag{$\mathbf{B}_m$}\] subject to the conditions ($\mathbf{A}_m$). We will show, over the course of Sections \ref{alg}--\ref{eq}, that ($\mathbf{B}_m$) does indeed hold for all $\bfx \in \Omega$ when $m$ is odd and $m \ge 3$. The corollary below is going to follow at once:

\begin{proposition} \label{odd} Suppose that $m$ is an odd positive integer with $m \ge 3$. Then the exponent of $2$ in the prime factorisation of the number of balanced grid colourings $B(m, n)$ is at least $2s_2(n)$. \end{proposition}

This result is of greatest relevance to the conjecture when, additionally, $s_2(m) = 2$; i.e., when $m$ is of the form $2^k + 1$. Though notice that Proposition \ref{odd} does not directly imply Theorem~\ref{2k1}; we must also make sure that $2s_2(n)$ is the exact exponent. On the other hand, our experience with the special case of Section~\ref{three} suggests that establishing the lower bound ought to be the most difficult part, and the rest of the proof might not require too much extra effort on top of that.

We figured out a quick derivation of $(\mathbf{B}_3)$ in Section \ref{three}. However, ($\mathbf{B}_m$) does not appear to be quite so straightforward when $m \ge 5$. In order to improve our understanding of the problem, we consider in Section \ref{alg} a wider class of analogously-shaped inequalities and we make some general observations about them.

\section{An Algorithmic Digression} \label{alg}

Our variables in this section will be $x_1$, $x_2$, $\ldots$, $x_h$, to be assigned nonnegative integer values.

For each $1 \le i \le p$, let $\bbD_i$ be a linear combination of the variables with nonnegative integer coefficients. We view these linear combinations as formal expressions over $x_1$, $x_2$, $\ldots$, $x_h$. Once each variable has been assigned a numerical value, we can evaluate all of $\bbD_1$, $\bbD_2$, $\ldots$, $\bbD_p$ as well.

We are interested in inequalities of the form \[c_1s_2(\bbD_1) + c_2s_2(\bbD_2) + \cdots + c_ps_2(\bbD_p) \ge 0 \tag{$\mathbf{I}$}\] with integer coefficients $c_1$, $c_2$, $\ldots$, $c_p$.

We may furthermore impose some conditions on the variables, of the form \[\bbE'_i = \bbE''_i \quad \text{for all $i$} \tag{$\mathbf{J}$}\] where $\bbE'_1$, $\bbE''_1$, $\ldots$, $\bbE'_q$, $\bbE''_q$ are linear combinations of $x_1$, $x_2$, $\ldots$, $x_h$ with nonnegative integer coefficients.

Given ($\mathbf{I}$) and ($\mathbf{J}$), we wish to determine whether ($\textbf{I}$) holds for all assignments of nonnegative integer values to the variables which satisfy ($\mathbf{J}$). If yes, then we want a proof; if not, a counterexample. Below, we will describe one algorithm which accomplishes this goal.

(Notice, though, that the purpose of this section is most of all motivational. Our eventual proof of Theorem \ref{2k1} will be ``properly mathematical'', intended to be comprehensible to a human mathematician without machine help.)

We begin with a quick refresher on how binary addition works. Suppose that we wish merely to add $x_1$, $x_2$, $\ldots$, $x_h$ together in binary. How do we do that, computationally?

Let $a$ be a nonnegative integer. We index the positions in the binary expansion of $a$ by $0$, $1$, $2$, $\ldots$, from right to left. So, if $a_i$ is the binary digit of $a$ at position $i$, then $a = \sum_i 2^ia_i$. We assume that the binary expansions of all numbers we are working with extend sufficiently far to the left -- which we ensure by padding them with meaningless leading zeroes as needed.

We can now compute the sum of $x_1$, $x_2$, $\ldots$, $x_h$ in binary by processing the positions one by one, in order, as follows: We keep track of the carry throughout. On each step $i$, let $r_i$ be the carry which the current position $i$ has received from the preceding position. Initially, we set $r_0 = 0$. To move forward, let $\sigma_i$ be the sum of the binary digits of $x_1$, $x_2$, $\ldots$, $x_h$ at position $i$. Then the binary digit at position $i$ in their sum will be \[(\sigma_i + r_i) \bmod 2;\] and the carry from the current position to the next one will be \[r_{i + 1} = \left\lfloor \frac{\sigma_i + r_i}{2} \right\rfloor.\]

The process stops when we exhaust all nonzero binary digits of $x_1$, $x_2$, $\ldots$, $x_h$ and the current carry vanishes. One simple but crucial observation is that, with $h$ summands, the carry can never exceed $h - 1$. This is straightforward to see by induction on the position.

Suppose now that, more ambitiously, we want to test whether a given assignment of nonnegative integer values to $x_1$, $x_2$, $\ldots$, $x_h$ satisfies ($\mathbf{J}$); and, if yes, then whether ($\mathbf{I}$) is true of it. Once again, we can do this by processing the positions one by one, in order. However, this time around we must keep track of a greater amount of information as we go. Specifically, we must keep track of the carries associated with all of the linear combinations $\bbD_1$, $\bbD_2$, $\ldots$, $\bbD_p$, $\bbE'_1$, $\bbE''_1$, $\ldots$,~$\bbE'_q$,~$\bbE''_q$.

Let $\delta_i$ be the sum of the coefficients of $\bbD_i$. Then the carry associated with $\bbD_i$ can never exceed $\delta_i - 1$. The carries associated with $\bbE'_i$ and $\bbE''_i$ are similarly bounded from above in terms of their own coefficient sums $\varepsilon'_i$ and $\varepsilon''_i$.

Given a positive integer $a$, let $[a]$ be the set $\{0, 1, \ldots, a - 1\}$. We define a \emph{state} to be any integer vector in $[\delta_1] \times [\delta_2] \times \cdots \times [\delta_p] \times [\varepsilon'_1] \times [\varepsilon''_1] \times \cdots \times [\varepsilon'_q] \times [\varepsilon''_q]$. The components of a state are meant to record the current carries associated with $\bbD_1$, $\bbD_2$, $\ldots$, $\bbD_p$, $\bbE'_1$, $\bbE''_1$, $\ldots$, $\bbE'_q$, $\bbE''_q$, respectively; i.e., a state packages together all of the information that we must keep track of. Given a state $\bfu$ and one of our linear combinations $\bbF$, we write $\bfu[\bbF]$ for the component of $\bfu$ associated with $\bbF$.

Our initial state is $\bfu_0 = \mathbf{0}$. On step $i$ of the computation, let $\bfu_i$ be the state formed by the carries which the current position $i$ has received from the preceding position. Let also $\bfalpha_i \in \{0, 1\}^h$ be the integer vector formed by the binary digits of $x_1$, $x_2$, $\ldots$, $x_h$ at position $i$, in order. We write $\bfalpha_i[x_j]$ for the $j$-th component of $\bfalpha_i$; i.e., the binary digit of $x_j$ at position $i$. We will record the newly computed binary digits at position $i$ of the linear combinations $\bbD_1$, $\bbD_2$, $\ldots$, $\bbD_p$, $\bbE'_1$, $\bbE''_1$, $\ldots$, $\bbE'_q$, $\bbE''_q$ into the integer vector $\bfbeta_i \in \{0, 1\}^{p + 2q}$. Its structure is similar to that of a state, and the meaning of $\bfbeta_i[\bbF]$ is analogous as well.

How do we move forward? Given a linear combination $\bbF$ of $x_1$, $x_2$, $\ldots$, $x_h$, we write $\bfalpha_i \llbr \bbF \rrbr$ for the corresponding linear combination of the components of $\bfalpha_i$; i.e., if $\bbF = f_1x_1 + f_2x_2 + \cdots + f_hx_h$, then $\bfalpha_i \llbr \bbF \rrbr = f_1\bfalpha_i[x_1] + f_2\bfalpha_i[x_2] + \cdots + f_h\bfalpha_i[x_h]$. The newly computed binary digits of our linear combinations are now given by \[\bfbeta_i[\bbF] = (\bfalpha_i \llbr \bbF \rrbr + \bfu_i[\bbF]) \bmod 2;\] and the carries from the current position to the next one are given by \[\bfu_{i + 1}[\bbF] = \left\lfloor \frac{\bfalpha_i \llbr \bbF \rrbr + \bfu_i[\bbF]}{2} \right\rfloor.\]

Once again, the process stops when we exhaust all nonzero binary digits of the variables and the current state returns to $\mathbf{0}$.

If at any step $i$ during the computation we observe that $\bfbeta_i[\bbE'_j] \neq \bfbeta_i[\bbE''_j]$ for some $j$, then we can tell immediately that this assignment of numerical values to $x_1$, $x_2$, $\ldots$, $x_h$ does not satisfy ($\mathbf{J}$) and must be discarded. Otherwise, the assignment does satisfy ($\mathbf{J}$), and our task becomes to test whether ($\mathbf{I}$) holds.

To this end, on each step $i$ we compute also a \emph{score} $w_i$, defined as \[w_i = c_1\bfbeta_i[\bbD_1] + c_2\bfbeta_i[\bbD_2] + \cdots + c_p\bfbeta_i[\bbD_p].\]

The sum of all scores is just the left-hand side of ($\mathbf{I}$). Hence, ($\mathbf{I}$) will hold for this particular assignment if and only if the sum of all scores is nonnegative.

Observe that, save for the final summation of the scores, the amount of information we must keep track of on each step of our computation is bounded -- there are only finitely many possible states. This is an indication that it might be helpful for our problem to construct some kind of finite automaton. We go on now to the details of the construction.

Let $\mathcal{G}^\dagger$ be the oriented multigraph whose vertices are all possible states and whose arcs are as specified below. We call $\mathbf{0}$ the \emph{root state} of $\mathcal{G}^\dagger$. For each state $\bfu$ and each integer vector $\bfalpha \in \{0, 1\}^h$, an arc $\gamma$ in $\mathcal{G}^\dagger$ points from $\bfu$ to the state $\bfv$ defined, as above, by \[\bfv[\bbF] = \left\lfloor \frac{\bfalpha \llbr \bbF \rrbr + \bfu[\bbF]}{2} \right\rfloor\] for all of our linear combinations $\bbF$. We call $\bfalpha$ the \emph{input} of $\gamma$. We also define the \emph{output} $\bfbeta$ of $\gamma$ by \[\bfbeta[\bbF] = (\bfalpha \llbr \bbF \rrbr + \bfu[\bbF]) \bmod 2\] for all of our linear combinations $\bbF$. Finally, the \emph{score} of $\gamma$ is defined by \[\score(\gamma) = c_1\bfbeta[\bbD_1] + c_2\bfbeta[\bbD_2] + \cdots + c_p\bfbeta[\bbD_p].\]

From $\mathcal{G}^\dagger$, we obtain $\mathcal{G}^\ddagger$ by removing all arcs $\gamma$ which do not satisfy the conditions $\bfbeta[\bbE'_i] = \bfbeta[\bbE''_i]$ for all $i$. Then, additionally, from $\mathcal{G}^\ddagger$ we obtain $\mathcal{G}$ by removing all states (as well as the arcs incident with them) which cannot be reached from the root, or for which the root cannot be reached from them.

We call $\mathcal{G}$ the \emph{carry graph} of ($\mathbf{I}$), subject to ($\mathbf{J}$). The preceding discussion makes it clear that ($\mathbf{I}$) holds, subject to ($\mathbf{J}$), if and only if every walk in $\mathcal{G}$ which begins and ends at the root is of nonnegative total score.

(Of course, in a practical implementation of the algorithm, the construction of the carry graph will be somewhat different. E.g., if we begin at the root and generate the states by levels, in order of their reachability from it along valid arcs, we can skip $\mathcal{G}^\dagger$ altogether and jump straight to a subgraph of $\mathcal{G}^\ddagger$. We can also get rid of the superfluous arcs in $\mathcal{G}$; i.e., if there are multiple arcs which point from $\bfu$ to $\bfv$, we can keep one of those with the lowest score and prune off the others.)

Since we have made sure that the root can reach all states of $\mathcal{G}$, and be reached by them, it is straightforward to see that a negative-score root-to-root walk will exist in $\mathcal{G}$ if and only if it contains a negative-score cycle. Given a concrete $\mathcal{G}$, we can easily test for this with the help of standard techniques from algorithmic graph theory. Below, we outline one such testing procedure.

We work by computing a series of arrays $W_0$, $W_1$, $W_2$, $\ldots$; for each $i$ and each state $\bfu$ of $\mathcal{G}$, the component $W_i[\bfu]$ of the $i$-th array will record the smallest possible score of a walk in $\mathcal{G}$ from the root to $\bfu$ with at most $i$ steps. Initially, $W_0[\mathbf{0}] = 0$ and $W_0[\bfu] = \infty$ for all non-root states $\bfu$. Given $W_i$, we compute the next array $W_{i + 1}$ as follows: For each state $\bfv$, let $\omega$ be the minimum of $W_i[\bfu] + \score(\gamma)$ over all states $\bfu$ and all arcs $\gamma$ which point from $\bfu$~to~$\bfv$. Then $W_{i + 1}[\bfv] = \min\{\omega, W_i[\bfv]\}$.

There are two possible outcomes of this process. If $\mathcal{G}$ contains a negative cycle, eventually we will obtain an array $W_i$ with $W_i[\mathbf{0}] < 0$. Then ($\mathbf{I}$) is false in general, subject to ($\mathbf{J}$). We can recover an explicit counterexample by stringing together the inputs of the arcs traversed by the corresponding negative-score root-to-root walk in $\mathcal{G}$. Otherwise, if $\mathcal{G}$ does not contain a negative cycle, at some point our arrays will stabilise. Or, in other words, eventually we will observe that $W_i = W_{i + 1}$ for some $i$. Then ($\mathbf{I}$) does indeed hold, subject to ($\mathbf{J}$).

This completes the description of our algorithm. It is quite entertaining to watch in action. Here, we limit ourselves to a few assorted examples. (The arc counts are all with the superfluous arcs pruned off.)

For the unconditional inequality $s_2(x) + s_2(y) \ge s_2(x + y)$, the carry graph has $2$ states and $4$ arcs, and the score arrays stabilise after $2$ iterations. The unconditional inequality $5s_2(5x) \ge s_2(x)$, which is false in general, yields a carry graph with $5$ states and $10$ arcs. Following $23$ iterations, the algorithm outputs the counterexample $x = 838861$; or $x = 11001100110011001101$ in binary. On the other hand, the carry graph of $(\textbf{B}_3)$ has $468$ states and $13296$ arcs, with the score arrays stabilising after $9$ iterations.

\section{Proof of the Inequality \texorpdfstring{($\mathbf{B}_m$)}{(B m)}} \label{ineq}

Suppose that $m$ is odd and $m \ge 3$. We now specialise the apparatus developed in Section \ref{alg} to the inequality ($\mathbf{B}_m$), subject to the system of conditions ($\mathbf{A}_m$). We omit the condition $\sumU = 2n$ of ($\mathbf{A}_m$) as it is non-homogeneous; but we are not going to need it anyway. We also omit those components of the states which correspond to the single-variable linear combinations in ($\mathbf{B}_m$), because the associated carries will always vanish. So each state of $\mathcal{G}$ consists of $3m + 1$ components which correspond to $\sumM_1$, $\ldots$, $\sumM_m$, $\sumU$, $\sumW_1$, $\sumB_1$, $\ldots$, $\sumW_m$, $\sumB_m$, in this order.

Let $\gamma$ be any arc in $\mathcal{G}$ which points from $\bfu$ to $\bfv$, with input $\bfalpha$ and output $\bfbeta$. Imagine for a moment that we can find some positive constant $c$ such that \[c \cdot \bfu[\sumU] + \score(\gamma) \ge c \cdot \bfv[\sumU]\] for all arcs $\gamma$ of $\mathcal{G}$. Then it would follow, for all states $\bfw$, that the score of a walk in $\mathcal{G}$ from the root to $\bfw$ can never fall below $c \cdot \bfw[\sumU]$. So, in particular, the score of a root-to-root walk can never become negative. By the discussion in Section \ref{alg}, this would guarantee our desired result.

We are not going to apply the above technique in its purest form. Instead, we will tweak it slightly. We define the \emph{surplus} of $\gamma$ by \[\surp(\gamma) = \score(\gamma) + 2 \cdot \bfu[\sumU] - 2 \cdot \bfv[\sumU].\]

Clearly, the total score of a closed walk in $\mathcal{G}$ equals its total surplus. The trouble is that the surplus is not always nonnegative. However, it is negative only very rarely -- and so it will be helpful to us anyway.

We call $\gamma$ \emph{positive} or \emph{negative} when its surplus is positive or negative, respectively. Otherwise, if the surplus vanishes, we call $\gamma$ \emph{neutral}. The key to the proof is that it is very difficult for an arc to turn out negative; rather, negativity occurs only in highly specific circumstances. For example, when $m = 3$, out of the $13296$ arcs in $\mathcal{G}$ a mere $15$ are negative. Our plan will be to show that, in a closed walk within $\mathcal{G}$, every negative arc must be preceded by a positive arc which balances it out.

First let us introduce some simplifying notation. Since $\bfu[\sumU]$ and $\bfalpha \llbr \sumU \rrbr$ are going to come up quite often, henceforth we denote them by $r$ and $\sigma$, respectively, as in the beginning of Section~\ref{alg}. We also let $r^\star = \bfv[\sumU]$; $d = \bfbeta[\sumU]$; $d_i = \bfbeta[\sumM_i]$ for all $i$; and $D = d_1 + d_2 + \cdots + d_m$.

We can now rewrite the surplus in a more convenient form. The score of $\gamma$, in our new notation, becomes $\sigma + D - 3d$. Furthermore, $r^\star = (\sigma + r - d)/2$ by the definition of $\mathcal{G}$. Hence, \begin{align*} \surp(\gamma) &= \sigma + D - 3d + 2r - \sigma - r + d\\ &= r + D - 2d. \end{align*}

We are going to need the following simple observation:

\begin{lemma} \label{wbm} For all $i$ and all states $\bfw$, it holds that $\bfw[\sumU] \ge \bfw[\sumW_i] + \bfw[\sumB_i] + \bfw[\sumM_i]$. \end{lemma}

\begin{myproof} Recall that all states are reachable from the root. We induct on the distance. The claim is obvious at the root itself. For the induction step, we assume that the claim is true of $\bfu$ and we must show that it is also true of $\bfv$. Since $\bfalpha \llbr \sumU \rrbr = \bfalpha \llbr \sumW_i \rrbr + \bfalpha \llbr \sumB_i \rrbr + \bfalpha \llbr \sumM_i \rrbr$, we get that \begin{align*} \bfv[\sumU] &= \left\lfloor \frac{\bfalpha \llbr \sumU \rrbr + \bfu[\sumU]}{2} \right\rfloor \ge \left\lfloor \frac{\bfalpha \llbr \sumW_i \rrbr + \bfalpha \llbr \sumB_i \rrbr + \bfalpha \llbr \sumM_i \rrbr + \bfu[\sumW_i] + \bfu[\sumB_i] + \bfu[\sumM_i]}{2} \right\rfloor\\ &\ge \left\lfloor \frac{\bfalpha \llbr \sumW_i \rrbr + \bfu[\sumW_i]}{2} \right\rfloor + \left\lfloor \frac{\bfalpha \llbr \sumB_i \rrbr + \bfu[\sumB_i]}{2} \right\rfloor + \left\lfloor \frac{\bfalpha \llbr \sumM_i \rrbr + \bfu[\sumM_i]}{2} \right\rfloor = \bfv[\sumW_i] + \bfv[\sumB_i] + \bfv[\sumM_i], \end{align*} and we are done. \end{myproof}

The necessary and sufficient condition for $\gamma$ to be negative is that $r + D < 2d$. But $d$ is a binary digit, meaning that $2d$ is always small. On the other hand, $r$ is a carry in the binary summation of a great number of nonnegative integers, and so we should expect it to often be quite large. In addition, $D$ is the sum of $m$ binary digits, and so once again it ought to be quite large usually. These considerations already indicate that the negative arcs are unlikely to be too common.

We go on now to pin down the structure of a negative arc more rigorously. Let $\bfw_\text{Neg}$ be the state of $\mathcal{G}$ where the carry associated with the sum of all variables is $1$ and the other carries vanish; i.e., $\bfw_\text{Neg}[\sumU] = 1$ and $\bfw_\text{Neg}[\bbF] = 0$ for all $i$ and all $\bbF$ out of $\sumW_i$, $\sumB_i$, $\sumM_i$.

\begin{lemma} \label{na} The surplus of a negative arc is always $-1$. Furthermore, all negative arcs emanate from the state $\bfw_\textnormal{Neg}$. \end{lemma} % negative arc

\begin{myproof} Suppose that $\gamma$ is negative. We consider the following cases, based on the value of $r$:

\smallskip

\emph{Case 1}. $r \ge 2$. Then $r + D \ge 2 \ge 2d$, a contradiction with $\gamma$ being negative.

\smallskip

\emph{Case 2}. $r = 0$. By Lemma \ref{wbm}, we get that $\bfu$ is the root. So $\bfbeta[\bbF] = \bfalpha \llbr \bbF \rrbr \bmod 2$ for all $i$ and all $\bbF$ out of $\sumW_i$, $\sumB_i$, $\sumM_i$, $\sumU$. However, $\bfbeta[\sumW_i] = \bfbeta[\sumB_i]$ by ($\mathbf{A}_m$). It follows that $\bfalpha \llbr \sumM_i \rrbr$ and $\bfalpha \llbr \sumU \rrbr = \bfalpha \llbr \sumW_i \rrbr + \bfalpha \llbr \sumB_i \rrbr + \bfalpha \llbr \sumM_i \rrbr$ are of the same parity. Hence, $d = \bfalpha \llbr \sumU \rrbr \bmod 2 = \bfalpha \llbr \sumM_i \rrbr \bmod 2 = d_i$ for all $i$; and so $D = md$. Since $m \ge 3$, this implies that $r + D \ge 2d$ -- once again in contradiction with $\gamma$ being negative.

\smallskip

\emph{Case 3}. $r = 1$. Then $d = 1$ and $D = 0$ because $\gamma$ is negative. The former makes $\sigma$ even as $d = (\sigma + r) \bmod 2$; while the latter implies that $d_1 = d_2 = \cdots = d_m = 0$.

Fix an index $i$. Since $\bfbeta[\sumW_i] = \bfbeta[\sumB_i]$ by ($\mathbf{A}_m$) and $d_i = 0$, we get that $\bfbeta[\sumW_i] + \bfbeta[\sumB_i] + \bfbeta[\sumM_i] = (\bfalpha \llbr \sumW_i \rrbr + \bfu[\sumW_i]) \bmod 2 + (\bfalpha \llbr \sumB_i \rrbr + \bfu[\sumB_i]) \bmod 2 + (\bfalpha \llbr \sumM_i \rrbr + \bfu[\sumM_i]) \bmod 2$ is even. Because $\sigma = \bfalpha \llbr \sumU \rrbr = \bfalpha \llbr \sumW_i \rrbr + \bfalpha \llbr \sumB_i \rrbr + \bfalpha \llbr \sumM_i \rrbr$ is even, too, we find that $\bfu[\sumW_i] + \bfu[\sumB_i] + \bfu[\sumM_i]$ must be even as well. However, $\bfu[\sumW_i] + \bfu[\sumB_i] + \bfu[\sumM_i] \le 1$ by Lemma \ref{wbm}. We conclude that, in fact, $\bfu[\sumW_i] = \bfu[\sumB_i] = \bfu[\sumM_i] = 0$ for all $i$. \end{myproof}

We are properly equipped by this point to carry out the rest of our plan.

\begin{lemma} \label{pa} Suppose that it is possible in $\mathcal{G}$ to continue from $\gamma$ along a negative arc. Then $\gamma$ itself must be positive. \end{lemma} % positive arc

\begin{myproof} Of course, $\bfv = \bfw_\text{Neg}$ and $r^\star = 1$ by Lemma \ref{na}. Since also $r^\star = \lfloor (\sigma + r)/2 \rfloor$, we get that \[\sigma \in \{2 - r, 3 - r\}. \tag{$\mathbf{P}$}\]

The necessary and sufficient condition for $\gamma$ to be positive is that $r + D > 2d$. We proceed to consider several cases based on the value of $r$, as follows:

\smallskip

\emph{Case 1}. $r \ge 3$. Then $r + D \ge 3 > 2d$, and so $\gamma$ is positive.

\smallskip

\emph{Case 2}. $r = 2$. If $d = 0$, then $\gamma$ is certainly positive. Suppose, otherwise, that $d = 1$. Since $d = (\sigma + r) \bmod 2$, we get that $\sigma$ is odd; and so $\sigma = 1$ by ($\mathbf{P}$).

Or, in other words, there is exactly one instance of the binary digit $1$ among the components of $\bfalpha$. Let $T$ be the unique type with $\bfalpha[x_T] = 1$. Since $m$ is odd, $T$ must necessarily contain the letter $\sfM$. (This is in fact the first step in the proof where we are making use of $m$ being odd.) Let $a$ be one index such that $T[a] = \sfM$; so $\bfalpha \llbr \sumM_a \rrbr = 1$.

On the other hand, $\bfv[\sumM_a] = 0$ because $\bfv = \bfw_\text{Neg}$; and also $\bfv[\sumM_a] = \lfloor (\bfalpha \llbr \sumM_a \rrbr + \bfu[\sumM_a])/2 \rfloor$. However, now $\bfu[\sumM_a] = 0$ and $d_a = (\bfalpha \llbr \sumM_a \rrbr + \bfu[\sumM_a]) \bmod 2 = 1$. Thus $D \ge d_a \ge 1$ and $r + D \ge 2 + 1 > 2d$, making $\gamma$ positive.

\smallskip

\emph{Case 3}. $r = 1$. Then $\gamma$ is certainly positive when $d = 0$. Suppose, otherwise, that $d = 1$. Since $d = (\sigma + r) \bmod 2$, we get that $\sigma$ is even; and so $\sigma = 2$ by ($\mathbf{P}$).

Or, in other words, there are exactly two instances of the binary digit $1$ among the components of $\bfalpha$. Let $S$ and $T$ be the two types with $\bfalpha[x_S] = \bfalpha[x_T] = 1$. Since $m$ is odd, both of them must contain the letter $\sfM$. Let $a$ and $b$ be two indices such that $S[a] = T[b] = \sfM$.

Suppose first that $a = b$. Then $\bfalpha \llbr \sumM_a \rrbr = 2$ and $\bfv[\sumM_a] = \lfloor (\bfalpha \llbr \sumM_a \rrbr + \bfu[\sumM_a])/2 \rfloor \ge 1$, in contradiction with $\bfv = \bfw_\text{Neg}$.

We conclude that necessarily $a \neq b$. Then, in order to avoid the same contradiction as in the previous paragraph, $\bfalpha$ must satisfy $\bfalpha \llbr \sumM_a \rrbr + \bfu[\sumM_a] = \bfalpha \llbr \sumM_b \rrbr + \bfu[\sumM_b] = 1$. However, now $d_a = (\bfalpha \llbr \sumM_a \rrbr + \bfu[\sumM_a]) \bmod 2 = 1$ and, similarly, $d_b = 1$ as well. Thus $D \ge d_a + d_b = 2$ and $r + D \ge 1 + 2 > 2d$, once again making $\gamma$ positive.

\smallskip

\emph{Case 4}. $r = 0$. Then, as in Case 2 of the proof of Lemma \ref{na}, we get that $\bfu$ is the root and $d_1 = \cdots = d_m = d$. Suppose, for the sake of contradiction, that $d_1 = \cdots = d_m = d = 0$. Since $\bfu$ is the root and $d_i = (\bfalpha \llbr \sumM_i \rrbr + \bfu[\sumM_i]) \bmod 2$, we get that $\bfalpha \llbr \sumM_i \rrbr$ is even for all $i$. Similarly, so is $\sigma$; and, consequently, $\sigma = 2$ by ($\mathbf{P}$).

Or, equivalently, there are exactly two instances of the binary digit $1$ among the components of $\bfalpha$. Hence, we may define the types $S$ and $T$ as well as the indices $a$ and $b$ just as in the preceding Case 3. Since $\bfalpha \llbr \sumM_a \rrbr$ is even and $\bfalpha \llbr \sumM_a \rrbr \in \{1, 2\}$, we find that $\bfalpha \llbr \sumM_a \rrbr = 2$. However, this entails the same contradiction with $\bfv = \bfw_\text{Neg}$ as in Case 3. We conclude that, in reality, $d_1 = \cdots = d_m = d = 1$. Thus $D = m$; $r + D = m \ge 3 > 2d$; and $\gamma$ becomes positive. \end{myproof}

What remains is to put the pieces together. Consider any closed walk in $\mathcal{G}$. (Though we only actually care about the ones which begin and end at the root.) By Lemma \ref{pa}, any negative arc in it is preceded by a positive one. By Lemma \ref{na}, the total surplus of these two arcs is nonnegative. Therefore, the total surplus over the entire closed walk -- which equals its total score -- must be nonnegative as well. By the discussion in Section \ref{alg}, this completes our proof of ($\mathbf{B}_m$), subject to ($\mathbf{A}_m$), when $m$ is odd and $m \ge 3$. As we observed in Section \ref{domino}, Proposition \ref{odd} follows immediately.

\section{The Case of \texorpdfstring{$m = 2^k + 1$}{m = 2 k + 1}} \label{eq}

We return now to the chain of reasoning which we began in Section \ref{domino}. By this point, we know that each $\bfx \in \Omega$ contributes a number of balanced colourings divisible by $2$ with exponent at least $2s_2(n)$. The rest of the proof of Theorem \ref{2k1} will not be too difficult in comparison. Since the case when $k = 1$ and $m = 3$ has been resolved in Section \ref{three} already, below we assume that $k \ge 2$ and $m \ge 5$.

We are not going to attempt a full description of the equality cases of ($\mathbf{B}_m$) because that would force us to do a lot of work with the carry graph $\mathcal{G}$ of Section \ref{ineq} whose structure is somewhat complicated. Instead, we limit ourselves to the following simple observation:

\begin{lemma} \label{ba} Consider any valid assignment of nonnegative integer values to the variables of~$(\mathbf{B}_m)$. Suppose that some nonzero integer occurs among these values $8$ or more times. Then $(\mathbf{B}_m)$ becomes strict. \end{lemma} % big arc

\begin{myproof} Consider the root-to-root walk in $\mathcal{G}$ which corresponds to the assignment. Since some nonzero value repeats $8$ or more times, at some point we must traverse an arc whose input contains among its components $8$ or more instances of the binary digit $1$. The state $\bfu$ that this arc points to will then satisfy $\bfu[\sumU] \ge 4$. Hence, for the next arc along which the walk continues from $\bfu$, the surplus is bounded from below by $2$. The material in Section \ref{ineq} now implies that the total surplus over the entire root-to-root walk will be strictly positive, which is equivalent to what we want. \end{myproof}

Call an $\bfx \in \Omega$ \emph{proper} when $x_T = 0$ for all types $T$ with $T[m] \neq \sfM$. This splits $\Omega$ into the subset $\Omega^\dagger$ of the proper $\bfx$ and the subset $\Omega^\ddagger$ of the non-proper ones. We proceed to show that the total contribution of $\Omega^\ddagger$ towards the number of balanced colourings is divisible by $2$ with exponent strictly greater than $2s_2(n)$. Our strategy for this will be to extract as much information as we can manage out of the automorphisms discussed in Section \ref{init}.

To begin with, if two distinct elements of $\Omega^\ddagger$ are interchanged by colour swapping, then we can discard the pair -- just as in Section \ref{three}. This leaves us with the subset $\Omega^\mathsection$ of those $\bfx \in \Omega^\ddagger$ which colour swapping preserves.

The next few ``batches'' we process will be defined in terms of column permutations. Let $\psi$ be any permutation of $\{1, 2, \ldots, m\}$. We say that two types $S$ and $T$ are related by $S = \psi(T)$ when $S[\psi(i)] = T[i]$ for all $i$. We also say that two multiplicity vectors $\bfx$ and $\bfy$ are related by $\bfx = \psi(\bfy)$ when $x_{\psi(T)} = y_T$ for all types $T$. (Here, $y_T$ is the multiplicity of $T$ as specified by~$\bfy$.) Of course, if $\bfx$ and $\bfy$ are related in this way, then they yield the same number of balanced colourings.

Construct a binary tree $\mathfrak{T}$ of height $k$ with $2^k$ leaves labelled $1$, $2$, $\ldots$, $2^k$. Each automorphism of $\mathfrak{T}$ induces a permutation of the leaves. Conversely, this permutation determines the automorphism uniquely. Clearly, all such permutations form a subgroup of size $2^{2^k - 1}$ within the symmetric group of degree $2^k$. Let $\Psi$ be the natural embedding of this subgroup within the symmetric group of degree $m$; i.e., each $\psi \in \Psi$ fixes $m$ and acts upon $1$, $2$, $\ldots$, $2^k$ in the same way as some automorphism of $\mathfrak{T}$ acts upon its leaves.

Partition $\Omega^\mathsection$ into orbits based on $\Psi$; i.e., two elements of $\Omega^\mathsection$ are in the same orbit if and only if they are related by an element of $\Psi$. Since the size of each orbit divides the size of $\Psi$, which is a power of $2$, we get that the size of a non-singleton orbit is always even. Hence, we may safely discard all such orbits. What remains is the subset $\Omega^\mathparagraph$ of those $\bfx \in \Omega^\mathsection$ for which the orbit is a singleton.

Partition the types, too, into orbits based on $\Psi$. Then any $\bfx \in \Omega^\mathparagraph$ must be constant over each such orbit; i.e., it must satisfy $x_S = x_T$ whenever $S$ and $T$ are related by an element of $\Psi$.

Fix an element $\bfx$ of $\Omega^\mathparagraph$. Since also $\bfx \in \Omega^\ddagger$, there exists some type $T$ with $T[m] \neq \sfM$ such that $x_T \neq 0$. By virtue of $T[m] \neq \sfM$, the letters $\sfW$ and $\sfB$ occur among the first $2^k$ letters of $T$ with opposite parities. It is straightforward to see, by induction on $k$, that this forces the orbit of $T$ to be of size at least $2^k$. (The bound is attained, e.g., when $T = \sfW\sfM\sfM \ldots \sfM\sfB$.)

The same reasoning applies also to the orbit of the colour-swapping complement $S$ of $T$. However, $\bfx$ assigns identical multiplicities to $S$ and $T$ as it is an element of $\Omega^\mathsection$. Furthermore, $\bfx$ must be constant over the orbits of $S$ and $T$ because it is an element of $\Omega^\mathparagraph$; and these two orbits must be distinct because $S[m] \neq T[m]$. We conclude that some nonzero value occurs among the components of $\bfx$ a minimum of $2^{k + 1} \ge 8$ times. (Recall that $k \ge 2$ and $m \ge 5$.) This allows us, by Lemma~\ref{ba}, to discard all of $\Omega^\mathparagraph$ as well.

We are free now to focus exclusively on $\Omega^\dagger$. Collectively, the proper multiplicity vectors correspond to those balanced colourings where all dominoes in the rightmost domino column are mixed. Split off this domino column; i.e., partition the original grid into two sub-grids $G'$ and $G''$ of sizes $2^{k + 1} \times 2n$ and $2 \times 2n$, respectively. Then each one of $G'$ and $G''$ is balanced individually. Thus the total number of balanced grid colourings associated with the proper multiplicity vectors becomes the product of $B(2^k, n)$ and $B(1, n)$. On the other hand, by Theorem \ref{2k}, the exponent of $2$ in each one of these two numbers equals $s_2(n)$. The proof of Theorem \ref{2k1} is complete.

\section{Further Work} \label{further}

We already know that the domino method certainly does not work off-the-shelf when $m$ is even or $s_2(m) \ge 3$. Though this does not preclude the possibility that some modification or augmentation of it might.

The smallest value of $m$ for which we have not sorted out the sequence $\mathcal{B}(m)$ is $m = 6$ -- which is still quite small. It would be nice to see proofs of the conjecture for other sequences of the form $\mathcal{B}(m)$, too. Especially significant would be any such proof with $s_2(m) \ge 3$. The smallest value of this kind is, of course, $m = 7$.

It might be worthwhile to look into higher-dimensional generalisations as well. Explicitly, fix a positive integer $d$ and consider the $d$-dimensional grid of size $2n_1 \times 2n_2 \times \cdots \times 2n_d$. The definition of a balanced colouring extends to it in an obvious manner. Does the exponent of $2$ in the number of $d$-dimensional balanced grid colourings $B(n_1, n_2, \ldots, n_d)$ exhibit any interesting~patterns?

\section*{Acknowledgements}

The author is thankful to Ankan Bhattacharya for telling him about his conjecture.

The present paper was written in the course of the author's PhD studies under the supervision of Professor Imre Leader. The author is thankful also to Prof.\ Leader for his unwavering support.

\end{document}